\documentclass{article}
\usepackage{amssymb,amsmath,amsthm,graphicx,cite}

\def\qed{\hfill {\hbox{${\vcenter{\vbox{               %HOLLOW SQUARE
   \hrule height 0.4pt\hbox{\vrule width 0.4pt height 6pt
   \kern5pt\vrule width 0.4pt}\hrule height 0.4pt}}}$}}}

\def\utr{\, \underline{\triangleright}\, }
\def\otr{\, \overline{\triangleright}\, }

\newtheorem{theorem}{Theorem}

\theoremstyle{definition}
\newtheorem{example}{Example}
\newtheorem{definition}{Definition}

\date{}

\title{\Large \textbf{Biquandle Brackets and Quivers}}

\author{Sam Nelson\footnote{Email: Sam.Nelson@cmc.edu. 
Partially supported by Simons Foundation collaboration grant 702597.}}

\begin{document}
\maketitle

\begin{abstract}
In this brief expository article we review the background for biquandle
bracket quivers -- including biquandles, biquandle homsets, biquandle coloring
quivers and biquandle brackets -- for a talk at the 70th Topology Symposium
at Nara Women's University in August 2023.
\end{abstract}

\parbox{6in} {\textsc{Keywords:} Quantum enhancements, biquandles, biquandle
counting invariants, biquandle brackets, trace diagrams, quivers

\smallskip

}

\section{\large\textbf{Introduction}}\label{I}

\textit{Biquandles} are algebraic structures with axioms motivated by the
Reidemeister moves in knot theory. Every oriented knot or link $L$ in 
$\mathbb{R}^2$ or $S^2$ has a \textit{fundamental biquandle} $\mathcal{B}(L)$
analogous to the fundamental group of a topological space. A finite biquandle 
$X$ determines an invariant of oriented knots and links called the 
\textit{biquandle homset invariant} $\mathrm{Hom}(\mathcal{B}(L),X)$
consisting of biquandle homomorphisms from the fundamental biquandle of $L$
to $X$. These homset elements have the property that they can be represented
visually as \textit{colorings} of a diagram of the oriented knot or link $L$
analogously to the way linear transformations between vector spaces are
represented by matrices. In particular, choosing a different diagram of
$L$ yields a different representation of the homset elements analogously to
the way choosing a different basis yields a different matrix representing 
the same linear transformation, with the role of change-of-basis matrices
played in the biquandle homset case by biquandle-colored Reidemeister moves.
See \cite{EN,K} for more.

From the homset, many useful computable invariants can be defined. The 
simplest of these is the cardinality of the homset, a non-negative 
integer-valued oriented link invariant known in the literature as the 
\textit{biquandle counting invariant}, denoted $\Phi_X^{\mathbb{Z}}(L)$.
Any invariant $\phi$ of biquandle-colored diagrams determines a
multiset-valued invariant of oriented knots and links called an 
\textit{enhancement} of the counting invariant. The first such examples use
a cohomology theory on the category of biquandles to define 2-\textit{cocycle}
and 3-\textit{cocycle} invariants, multiset-valued invariants which can be 
transformed into polynomial invariants which evaluate to $\Phi_X^{\mathbb{Z}}(L)$
at zero but in general are stronger invariants; see \cite{CJKLS}. 
A more recent example of an 
enhancement is the \textit{biquandle coloring quiver} associated to a 
subset $S$ of the automorphism group of the finite biquandle $X$; such a subset
determines a directed graph-valued enhancement of the counting invariant.
Since directed graphs (also known as \textit{quivers}) are categories, this
enhancement is actually a \textit{categorification} of the counting invariant.
Several new polynomial invariants can be obtain from these quivers via
different forms of decategorification; see \cite{CN,CN2,CN3,IN,FN}

Our main interest in this manuscript is in a family of enhancements
known as \textit{biquandle brackets}. These are a type of \textit{quantum 
enhancement}, i.e. quantum invariants of biquandle-colored oriented knots and 
links representing elements of the homset $\mathrm{Hom}(\mathcal{B}(L),X)$. 
In this manuscript and its associated talk, we review (in a fairly 
self-contained way) biquandles and biquandle quivers, leading up to biquandle 
brackets and biquandle bracket quivers. In Section \ref{B}
we review the basics of biquandles and the biquandle homset invariant. In 
Section \ref{Q} we review the quiver categorification of the biquandle homset 
invariant. In Section \ref{BB} we review biquandle brackets and their 
categorification via quivers. We conclude in Section \ref{Q} with a few words 
about current work in this area.

\section{\large\textbf{Knots, Biquandles and Biquandle Coloring}}\label{B}

A \textit{knot} is an embedding $K:S^1\to \mathbb{R}^3$ or $K:S^1\to S^3$;
we may also consider the image of such an embedding as a knot. A \textit{link}
is a disjoint union of knots. A knot is \textit{tame} if the embedding is 
piecewise-linear or (equivalently) $C^{\infty}$.
In the 1920s, Kurt Reidemeister, a mathematician at the Georg August 
Universit\"at in G\"ottingen, Germany, proved that two diagrams 
represent ambient isotopic tame knots in $S^3$ iff they are related by a
sequence of the following moves, now known as \textit{Reidemeister moves}:
\[\includegraphics{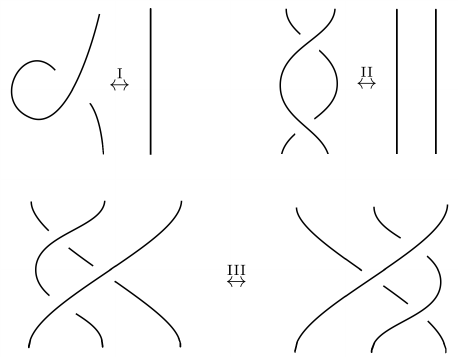}.\]
For the remainder of this Manuscript, we will consider only tame knots and links.
See \cite{R} for more.

Thus, to prove that two diagrams represent the same knot, it suffices to
identify a sequence of Reidemeister moves changing one diagram to the other.
To prove that two diagrams represent different knots, we can identify
\textit{knot invariants}, functions $f:\mathcal{D}\to X$ taking knot diagrams 
as inputs and giving output values in some set $X$ (which could be integers, 
polynomials, graphs etc.) with the property that diagrams differing by
Reidemeister moves have the same function value, i.e.,
\[D\sim D'\Rightarrow f(D)=f(D').\]
Then if two knot diagrams share the same value of an invariant, it says 
nothing -- the two diagrams could represent the same knot or they could just 
coincidentally have the same value; however, if two diagrams have different 
values of an invariant, then they cannot be related by Reidemeister moves 
and hence must represent different knots.

To define invariants, we need to find things that are not changed by 
Reidemeister moves. One way to do this is to use the power of universal 
algebra, creating algebraic structures with axioms derived from the 
Reidemeister moves. The example of interest for this manuscript is the 
algebraic structure known as \textit{biquandles} (see \cite{EN,K} for more).

\begin{definition} 
A \textit{biquandle} is a set $X$ with a pair of binary operations
$\utr,\otr:X\times X\to X$ satisfying the axioms
\begin{itemize}
\item[(i)] For all $x\in X$, $x\utr x=x\otr x$,
\item[(ii)] For all $y$ in $X$, the maps $\alpha_y,\beta_y:X\to X$ defined by 
$\alpha_y(x)=x\otr y$ and $\beta_y(x)=x\utr y$ and the map 
$S:X\times X\to X\times X$ defined by $S(x,y)=(y\otr x, x\utr y)$ are invertible, and
\item[(iii)] For all $x,y,z\in X$, we have
\[\begin{array}{rcl}
(x\utr y)\utr (z\utr y) & = & (x\utr z)\utr (y\otr z)  \\
(x\utr y)\otr (z\utr y) & = & (x\otr z)\utr (y\otr z)  \\
(x\otr y)\otr (z\otr y) & = & (x\otr z)\otr (y\utr z)
\end{array}.\]
\end{itemize}
\end{definition}

These axioms result from interpreting the elements of $X$ as labels or
``colors'' for the semiarcs (segments between crossings) in an oriented 
knot diagram and interpreting the operations as over- and under-crossing 
rules as shown:
\[\scalebox{2}{\includegraphics{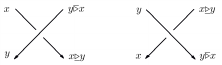}}\]

Then axiom (i) expresses the four oriented versions of the Reidemeister I move:
\[\includegraphics{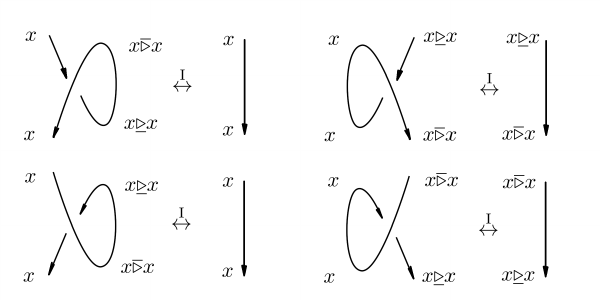}.\]
Similarly axiom (ii) expresses the four oriented Reidemeister II moves
and axiom (iii) expresses the all-positive Reidemeister III move.

Since these moves form a generating set of oriented Reidemeister moves (see 
\cite{P}), we have:

\begin{theorem}
If $X$ is a biquandle and $D$ and $D'$ are oriented knot or link diagrams related
by Reidemeister moves, there is a one-to-one correspondence between the
set of $X$-colorings of $D$ and the set of $X$-colorings of $D'$.
\end{theorem}

See \cite{EN,K} for more details.

\begin{definition}
A map $f:X\to Y$ between biquandles is a \textit{homomorphism} if for all
$x,y\in X$ we have
\[f(x\utr y)=f(x)\utr f(y)\ \mathrm{and} \ f(x\otr y)=f(x)\otr f(y).\]
A bijective homomorphism is an \textit{isomorphism}.
\end{definition}

\begin{example}
Given a diagram representing an oriented knot or link diagram $L$, we 
define the \textit{fundamental biquandle} $\mathcal{B}(L)$ in the following way:
\begin{itemize}
\item We form a set of \textit{generators} consisting of a symbol for each
semiarc in the diagram, 
\item We form a set of \textit{biquandle words} including generators and 
expressions of the forms $u\utr v$, $u\otr v$, $\alpha_u^{-1}(v)$, 
$\beta_u^{-1}(v)$ and $S_i^{-1}(x_j,x_k)$ for $i\in \{1,2\}$ (we interpret 
$S_i^{-1}(u,v)$ as the $i$th component of $S^{-1}(u,v)$) where $u,v$ are 
biquandle words,
\item Then $\mathcal{B}(L)$ is the set of equivalence classes of biquandle
words under the equivalence relation generated by the biquandle axioms and the
crossing relations.
\end{itemize}
The set of generators and crossing relations (taking the biquandle axioms
as understood) is known as a \textit{biquandle presentation}.
Reidemeister moves on diagrams produce Tietze moves on the 
presentation, resulting in isomorphic biquandles; hence, the fundamental 
biquandle does not depend on our choice of diagram for $L$ and is a link 
invariant.
\end{example}

An $X$-coloring of an oriented knot or link diagram $D$ is an assignment 
of an image in $X$ to each generator of $\mathcal{B}(L)$ satisfying
the crossing relations and hence determines a unique biquandle homomorphism
$f:\mathcal{B}(K)\to X$. Thus the set of biquandle colorings of a diagram $D$
of an oriented knot or link $L$ can be identified with the homset
$\mathrm{Hom}(\mathcal{B}(L),X)$.
In particular, a choice of diagram for $L$ is analogous to a choice of
basis for a vector space -- just as a linear transformation $f:X\to Y$ is 
determined by a choice of output vector for each input basis vector, a 
biquandle homset element is determined by a choice of output color in $X$
for each semiarc in our chosen diagram of $L$ such that the crossing relations
are satisfied in $X$. Different choices of diagram
for $L$ yield different representations of the same homset element just
as different choices of basis yield different matrices representing the
same linear transformation, with the role of change-of-basis matrices 
played by $X$-colored Reidemeister moves.

\begin{example} \label{ex1}
Let $X$ be the biquandle $\mathbb{Z}_2=\{0,1\}$ with $x\utr y=x\otr y=x+1$ 
and consider the Hopf link $L$. The homset $\mathrm{Hom}(\mathcal{B}(L),X)$
has four elements which we can represent with the set of $X$-colored diagrams
\[\includegraphics{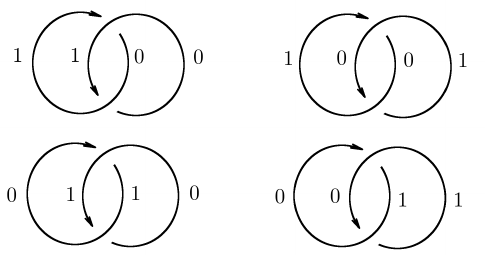}.\]
\end{example}

\begin{example}
Each of the diagrams below represents the same homset element in 
$\mathrm{Hom}(\mathcal{B}(L),X)$ where $L$ is the unknot and $X$ is the 
biquandle in Example \ref{ex1}.
\[\includegraphics{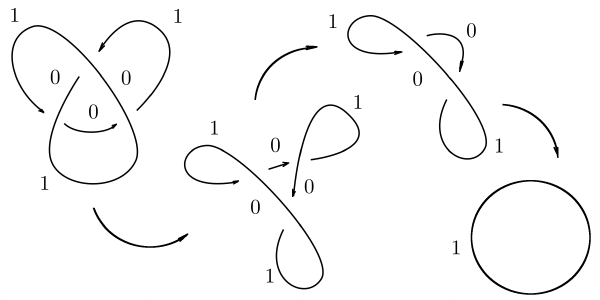}\]
\end{example}

See \cite{EN} for more.

\section{\large\textbf{Quivers and Categorification}}\label{Q}

\textit{Categories} are algebraic structures consisting of a collection
of \textit{objects} and for each pair of objects $X,Y$ a set of 
\textit{morphisms} denoted $\mathrm{Hom}(X,Y)=\{f:X\to Y\}$ satisfying certain
axioms. Examples include the category of sets and functions, the category
of vector spaces and linear transformations, the category of groups and
group homomorphisms, the category of topological spaces and continuous maps, 
etc. Indeed, many subject areas of mathematics can be described in 
category-theoretic terms. Connections between these subject areas can often
be formalized as maps known as \textit{functors} between categories, e.g. 
the fundamental group functor $\pi_1$ which transforms topological questions 
into group-theoretic questions. Similarly, functors known as \textit{homology 
theories} measure information about topological or algebraic structures in 
terms of chain complexes and homology groups.

In the late 1990s, mathematical physicists such as Louis Crane and (my former
colleague) John Baez proposed a program of \textit{categorification} wherein
simpler structures are replaced with richer categorical structures. For
example, we might replace natural numbers $n$ with vector spaces of 
dimension $n$, addition with direct sum and multiplication 
with tensor product. As a result, we get a richer and more powerful structure --
while integers can be either equal or not, vector spaces may be identical, 
distinct but isomorphic, or different. From these richer structures we can
recover the original simpler structures and potentially get new ones via 
\textit{decategorification}, e.g. taking the dimension of the vector spaces.
Famous examples in recent decades include Khovanov homology (categorifying 
the Jones polynomial) and Knot Floer homology (categorifying the Alexander 
polynomial). 

Directed graphs, sometimes known as \textit{quivers} (as collections of 
arrows), form categories with vertices as objects and directed paths as 
morphisms. In recent years my students, professional collaborators and
I have used quivers to categorify several enhancements of the biquandle
counting invariant; see \cite{CN,CN2,CN3,IN,KNS} for more.

For this section we will focus on one example, the 
biquandle coloring quiver. Let $X$ be a finite biquandle and $D$ an 
oriented knot or link diagram with a choice of coloring $X$-coloring, and 
let $f:X\to X$ be a biquandle endomorphism. Then applying $f$ to each color
on the diagram $D$ yields another assignment of elements of $X$ to generators
of $\mathcal{B}(D)$, and the fact that $f$ is a quandle endomorphism implies
that this assignment is a valid $X$-coloring. 

\[\scalebox{2}{\includegraphics{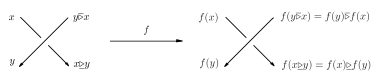}}\]

In particular, each endomorphism $f$ determines an arrow from one coloring
in the homset to another. Thus, a set $S$ of endomorphisms of $X$ determines
a quiver structure on the homset which we call the \textit{biquandle coloring
quiver} associated to $S$. Changing the diagram $D$ by $X$-colored 
Reidemeister moves does not change the homset or the resulting quiver, and 
hence the quiver is an invariant of knots and links for each $S$. In 
particular the case of $S=\emptyset$ can be identified with the original 
homset invariant. The case of $S=\mathrm{Hom}(X,X)$, the entire set of 
endomorphisms, is called the \textit{full quiver}. 

\begin{example}
In our previous Example \ref{ex1}, the identity map and the map $f(x)=x+1$ 
are the only endomorphisms. Then the full biquandle coloring quiver looks like
\[\includegraphics{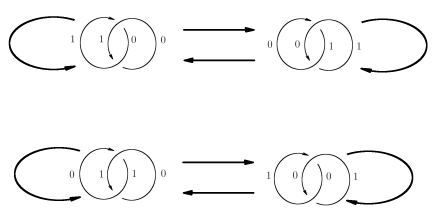}.\]
\end{example}

Such a quiver is a category and hence the biquandle coloring quiver is a 
categorification of the biquandle homset invariant; it decategorifies to the 
original homset by deleting the arrows. Other polynomial invariants such 
as the \textit{in-degree polynomial} can be obtained as decategorifcations; 
see \cite{CN,CN2} for more.

\section{\large\textbf{Biquandle Brackets}}\label{BB}

In the 1980s, Vaughn Jones won the Fields medal for his introduction of 
a powerful knot invariant known as the \textit{Jones polynomial}. 
The invariant can be defined using the \textit{Kauffman bracket skein relation}
which we may understand as thinking of crossings as linear combinations
of \textit{smoothings} or diagrams with crossings removed. More precisely,
applying the skein relation to a diagram with a crossing replaces the
diagram with a linear combination of two diagrams with the crossing removed
and the resulting endpoints connected in the two possible ways with 
coefficients in $\mathbb{Z}[A^{\pm 1}]$. 
Recursively applying this procedure to all crossings yields a linear 
combination of diagrams without crossings, known as \textit{Kauffman states}.
Evaluating each Kauffman state as $\delta^{k-1}$ where $k$ is the number of 
components and $\delta=-A^2-A^{-2}$ yields the \textit{unnormalized Kauffman
bracket polynomial}, and multiplying by $(-A^3)^{-wr}$ where $wr$ is the 
\textit{writhe} or number of positive crossings minus the number of negative
crossings yields the Kauffman bracket polynomial, which can be easily checked 
to be invariant under Reidemeister moves.

For biquandle brackets we want to generalize the Kauffman bracket to
the case of biquandle colored knots and links. To do this, we replace the
Kauffman bracket skein relation with the skein relations
\[\includegraphics{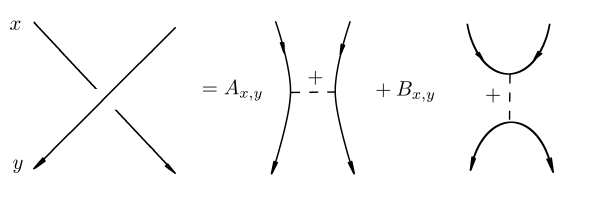}\]
\[\includegraphics{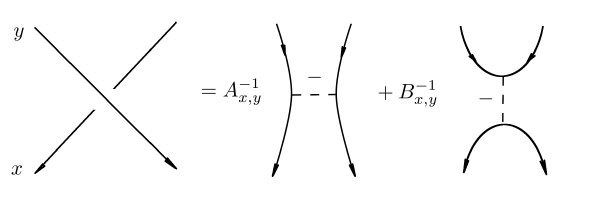}\]
yielding \textit{trace diagrams} with dashed signed edges called a 
\textit{traces}. Traces act like crossings, enabling us to retain the 
biquandle colors in a homset element diagram. We now have skein coefficients
$A_{x,y}$ and $B_{x,y}$ which are functions of the biquandle colors $x,y\in X$
on the left side of the crossing. Deleting the traces and colors in a 
fully-smoothed diagram yields the usual Kauffman states, and we require
the elements $-A_{x,y}B_{x,y}^{-1}-A_{x,y}^{-1}B_{x,y}$ to be equal, with their
common value denoted as $\delta$.

The Reidemeister moves then impose some conditions on the skein coefficients:
\[\begin{array}{rcl}
A_{x,y}A_{y,z}A_{x\utr y,z\otr y} & = & A_{x,z}A_{y\otr x,z\otr x}A_{x\utr z,y\utr z} \\
A_{x,y}B_{y,z}B_{x\utr y,z\otr y} & = & B_{x,z}B_{y\otr x,z\otr x}A_{x\utr z,y\utr z} \\
B_{x,y}A_{y,z}B_{x\utr y,z\otr y} & = & B_{x,z}A_{y\otr x,z\otr x}B_{x\utr z,y\utr z} \\
A_{x,y}A_{y,z}B_{x\utr y,z\otr y} & = & 
A_{x,z}B_{y\otr x,z\otr x}A_{x\utr z,y\utr z} 
+A_{x,z}A_{y\otr x,z\otr x}B_{x\utr z,y\utr z} \\ 
& & +\delta A_{x,z}B_{y\otr x,z\otr x}B_{x\utr z,y\utr z} 
+B_{x,z}B_{y\otr x,z\otr x}B_{x\utr z,y\utr z} \\
B_{x,y}A_{y,z}A_{x\utr y,z\otr y} 
+A_{x,y}B_{y,z}A_{x\utr y,z\otr y} & & \\
+\delta B_{x,y}B_{y,z}A_{x\utr y,z\otr y} 
+B_{x,y}B_{y,z}B_{x\utr y,z\otr y}  
& = & B_{x,z}A_{y\otr x,z\otr x}A_{x\utr z,y\utr z}. \\
\end{array}\]

A biquandle bracket $\beta$ defined on a biquandle $X$ and a commutative 
unital ring $R$ then consists of a pair of maps $A,B:X\times X\to R^{\times}$
such that the above equations are satisfied and the values 
$\delta=-A_{x,y}B_{x,y}^{-1}-A_{x,y}^{-1}B_{x,y}$ and $w=-A_{x,x}^2B_{x,x}^{-1}$
are the same for all $x,y\in X$. Then for each $X$-coloring of an oriented 
link diagram $D$, the sum of the products of state values times smoothing 
coefficients times writhe adjustment $w^{-wr}$ over the set of all Kauffman 
states, known as the \textit{state-sum} value, is unchanged by $X$-colored 
Reidemeister moves, and the multiset of such values over the complete homset
is an invariant of oriented knots and links.

\begin{example}
Let $X$ be the biquandle from our previous Example \ref{ex1} and let 
$R=\mathbb{Z}_5$. Then one can verify that the coefficient tables
\[
\begin{array}{r|rr}
A & 1 & 0 \\ \hline
1 & 1 & 3 \\
0 & 4 & 1 
\end{array}\quad
\begin{array}{r|rr}
B & 1 & 0 \\ \hline
1 & 4 & 2 \\
0 & 1 & 4
\end{array}\quad
\]
satisfy the biquandle bracket axioms with $\delta=2$ and $w=1$. 
Then the Hopf link with its four 
$X$-colorings has biquandle bracket multiset $\{3,3,4,4\}$:
\[\includegraphics{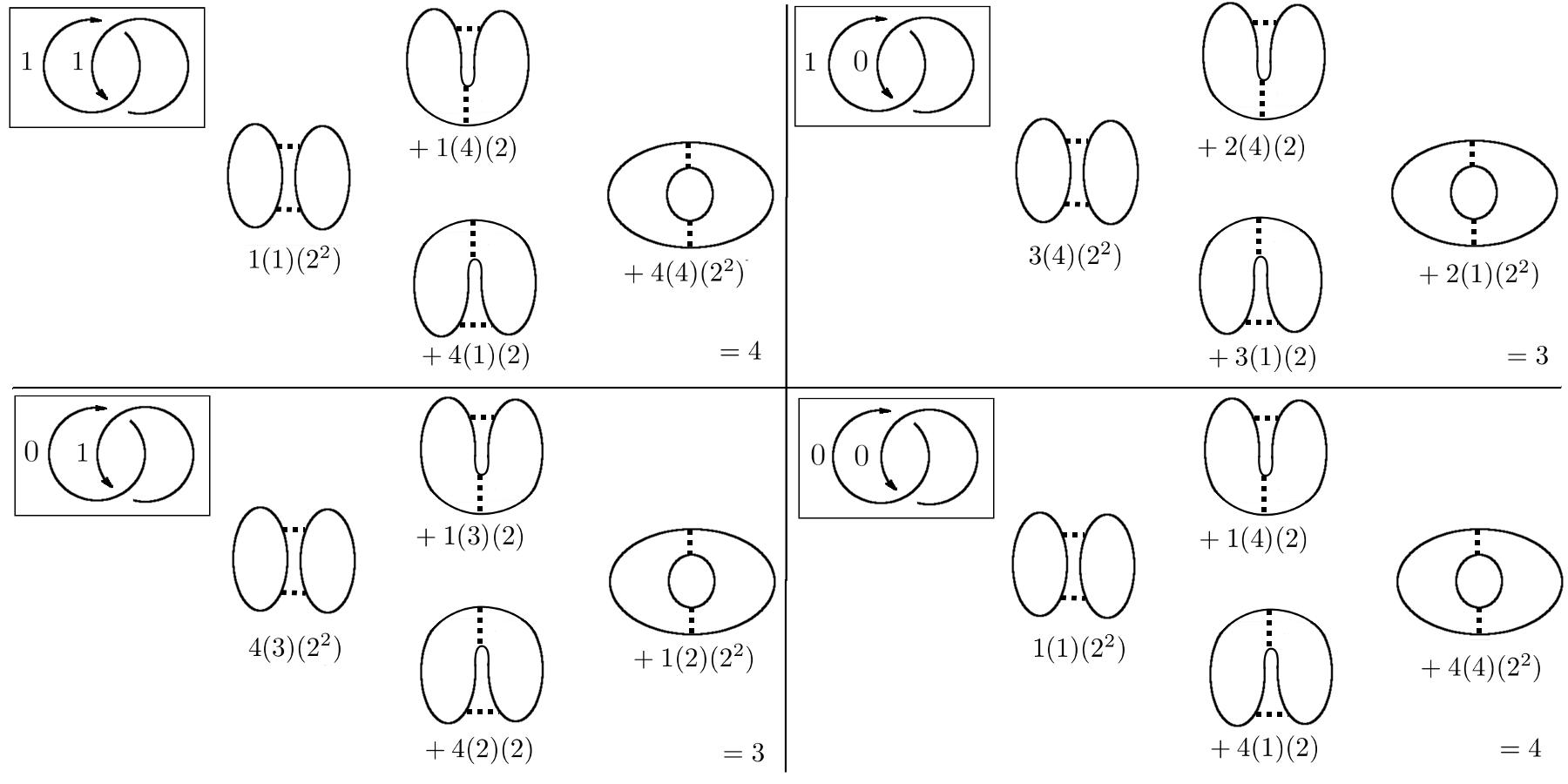}\]
\end{example}

See \cite{NO,NOR} for more details and examples.

Finally, we can categorify the biquandle bracket invariant using the 
biquandle coloring quiver. More precisely, where the biquandle bracket
invariant is a multiset of $\beta$-values over the biquandle homset,
a choice of set of biquandle endomorphisms $S$ gives us an invariant 
quiver with vertices weighted with $\beta$-values which we call a 
\textit{biquandle bracket quiver}. 

\begin{example}\label{ex2}
Continuing with the Hopf link and biquandle bracket from the previous 
examples, we have the following biquandle bracket quiver:
\[\includegraphics{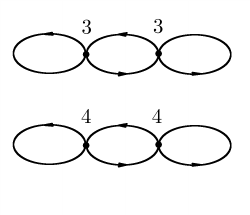}\]
\end{example}

These quivers can become very large and complex very quickly, so for useful 
invariants it is helpful to decategorify in various ways to obtain
easy-to-use polynomial invariants.

\begin{example}
Given a biquandle bracket quiver, we can sum over the set of arrows
terms of the form $s^{S(a)}t^{T(a)}$ where $S(a)$ and $T(a)$ are the weights
at the source and target of the arrow respectively. In the biquandle bracket
quiver in Example \ref{ex2} above, we get $4s^3t^3+4s^4t^4$ for this
polynomial decategorification.
\end{example}

\begin{example}
Another decategorification uses the fact that while the out-degree of 
every vertex in a biquandle bracket quiver is the same, the in-degrees can be 
(and generally are) different. Thus, we can sum over the set of vertices
terms of the form $u^{\beta(v)}w^{\mathrm{deg}_+(v)}$ where $\beta(v)$ is the
weight at the vertex $v$ and $\mathrm{deg}_+(v)$ is the in-degree of the 
vertex $v$. Then in our Example \ref{ex2} we obtain $2u^3w^2+2u^4w^2$ for this
decategorification.
\end{example}

See \cite{FN} for more.

\section{\large\textbf{Current and Future Work}} \label{Q}

I am currently working with several collaborators around the world on 
generalizations and extensions of these biquandle bracket quiver invariants,
and other groups are also working on biquandle brackets. In \cite{HVW} it is
shown that many biquandle brackets are cohomologous to the Jones polynomial,
meaning that biquandle bracket quivers yield a family of new categorifications
of the Jones polynomial quite unlike Khovanov homology it is various forms.
In \cite{IM} biquandle brackets are generalized to picture-valued invariants.

Faster algorithms for finding biquandle brackets are of great interest, as 
are new methods of finding biquandle brackets over infinite coefficient 
rings. Work currently in preparation (check \texttt{arXiv.org} soon!)
includes \textit{biquandle power brackets}, a major advancement over
standard biquandle brackets in which the value $\delta$ of a Kauffman state 
component can vary as a function of the biquandle colors it contains, as well
as \textit{biquandle-colored Conway algebras} and extensions of the biquandle 
bracket to other algebraic structures.

\bibliography{bbnq}{}

\begin{thebibliography}{10}

\bibitem{CJKLS}
J.~S. Carter, D.~Jelsovsky, S.~Kamada, L.~Langford, and M.~Saito.
\newblock State-sum invariants of knotted curves and surfaces from quandle
  cohomology.
\newblock {\em Electron. Res. Announc. Amer. Math. Soc.}, 5:146--156, 1999.

\bibitem{CN}
J.~Ceniceros and S.~Nelson.
\newblock Psyquandle coloring quivers.
\newblock {\em Preprint, arXiv:2107.05668}.

\bibitem{CN3}
K.~Cho and S.~Nelson.
\newblock Quandle cocycle quivers.
\newblock {\em Topology Appl.}, 268:106908, 10, 2019.

\bibitem{CN2}
K.~Cho and S.~Nelson.
\newblock Quandle coloring quivers.
\newblock {\em J. Knot Theory Ramifications}, 28(1):1950001, 12, 2019.

\bibitem{EN}
M.~Elhamdadi and S.~Nelson.
\newblock {\em Quandles---an introduction to the algebra of knots}, volume~74
  of {\em Student Mathematical Library}.
\newblock American Mathematical Society, Providence, RI, 2015.

\bibitem{FN}
P.~C. Falkenburg and S.~Nelson.
\newblock Biquandle bracket quivers.
\newblock {\em J. Knot Theory Ramifications (to appear)}.

\bibitem{HVW}
W.~Hoffer, A.~Vengal, and V.~Winstein.
\newblock The structure of biquandle brackets.
\newblock {\em J. Knot Theory Ramifications}, 29(6):2050042, 13, 2020.

\bibitem{IM}
D.~P. Ilyutko and V.~O. Manturov.
\newblock Picture-valued parity-biquandle bracket.
\newblock {\em J. Knot Theory Ramifications}, 29(2):2040004, 22, 2020.

\bibitem{IN}
K.~Istanbouli and S.~Nelson.
\newblock Quandle module quivers.
\newblock {\em J. Knot Theory Ramifications}, 29(12):2050084, 14, 2020.

\bibitem{K}
L.~H. Kauffman and D.~Radford.
\newblock Bi-oriented quantum algebras, and a generalized {A}lexander
  polynomial for virtual links.
\newblock In {\em Diagrammatic morphisms and applications ({S}an {F}rancisco,
  {CA}, 2000)}, volume 318 of {\em Contemp. Math.}, pages 113--140. Amer. Math.
  Soc., Providence, RI, 2003.

\bibitem{KNS}
J.~Kim, S.~Nelson, and M.~Seo.
\newblock Quandle coloring quivers of surface-links.
\newblock {\em J. Knot Theory Ramifications}, 30(1):Paper No. 2150002, 13,
  2021.

\bibitem{NOR}
S.~Nelson, M.~E. Orrison, and V.~Rivera.
\newblock Quantum enhancements and biquandle brackets.
\newblock {\em J. Knot Theory Ramifications}, 26(5):1750034, 24, 2017.

\bibitem{NO}
S.~Nelson and N.~Oyamaguchi.
\newblock Trace diagrams and biquandle brackets.
\newblock {\em Internat. J. Math.}, 28(14):1750104, 24, 2017.

\bibitem{P}
M.~Polyak.
\newblock Minimal generating sets of {R}eidemeister moves.
\newblock {\em Quantum Topol.}, 1(4):399--411, 2010.

\bibitem{R}
K.~Reidemeister.
\newblock {\em Knotentheorie}.
\newblock Springer-Verlag, Berlin-New York, 1974.
\newblock Reprint.

\end{thebibliography}
\bibliographystyle{abbrv}

\bigskip

\noindent
\textsc{Department of Mathematical Sciences \\
Claremont McKenna College \\
850 Columbia Ave. \\
Claremont, CA 91711}

\end{document}